\documentclass[12pt]{article}
\usepackage{setspace}

\usepackage{amsmath}
\usepackage{framed}
\def\qed{\rule[0mm]{.6em}{.6em}}

\begin{document}

\setcounter{secnumdepth}{3}
\setcounter{tocdepth}{3}
\newtheorem{Lemma}{Lemma}[section]
\newtheorem{Theorem}{Theorem}[section]
\newtheorem{Corollary}{Corollary}[section]
\newenvironment{Proof}{{\bf Proof }}{\qed}

\begin{center}

\section*{G\"{o}del for Goldilocks: A Rigorous, Streamlined Proof of (a variant of) 
G\"{o}del's First Incompleteness Theorem\footnote{This exposition requires minimal background. Other than
common things, you only need to know
what an integer is; what a function is; and what a computer program is. 
You do need several hours, and you need to focus---the material is concrete and understandable, but it is not trivial. 
This material was the basis for the first two lectures 
of my course offering ``The Theory of Computation" (a sophomore/junior level course) in October 2014.}}

{\bf Dan Gusfield} 

Department of Computer Science, UC Davis

August 2014,
revised November 15, 2014
\end{center}

\medskip

\section{Introduction: Why I wrote this} 

G\"{o}del's famous incompleteness theorems (there are two of them) 
concern the ability of a formal system to state and derive
all true statements, and only true statements, in some fixed domain; and concern the ability of logic to determine if a 
formal system has that property. They were developed in the early 1930s. 
Very loosely, the first theorem says that in any ``sufficiently rich" formal proof system where
it is not possible to prove a false statement about {\it arithmetic}, there will also be true statements about arithmetic that cannot be proved.

Most discussions of G\"{o}del's theorems fall into one of two types: either they emphasize perceived cultural and philosophical meanings of the
theorems, and perhaps sketch some of the ideas of the proofs, usually relating G\"{o}del's proofs to riddles and paradoxes, 
but do not attempt rigorous, complete proofs;
or they do present rigorous proofs, 
but in the traditional style of mathematical logic, with all of its heavy notation,
difficult definitions, technical issues in G\"{o}del's original approach, and connections
to broader logical theory before and after G\"{o}del.
Many people are frustrated
by these two extreme types of expositions\footnote{To verify this, just randomly search the web
for questions about G\"{o}del's theorem.} and want a short, straight-forward, rigorous proof that they can understand.

Over time, various people have realized that somewhat weaker, but still meaningful, variants of G\"{o}del's first incompleteness theorem can 
be rigorously proved by simpler arguments based on notions of computability. This approach avoids 
the heavy machinery of mathematical logic at one extreme; and does not rely on analogies, paradoxes, 
philosophical discussions or hand-waiving, at the other extreme. 
This is the just-right {\it Goldilocks} approach.
However, the available expositions of this middle approach have still been aimed at a relatively sophisticated audience, and 
have either been too brief,\footnote{For example, in Scott Aaronson's book {\it Quantum Computing Since Democritus}.}
or have been embedded in larger, more involved discussions.\footnote{For example, Sipser's excellent
book on the Theory of Computation, where the exposition of G\"{o}del's theorem 
relies on an understanding of Turing machines and the Undecidability of the Halting problem. Another example is {\it An Introduction to G\"{o}del's
Theorems} by Peter Smith, which develops much more logical machinary before proving a variant of G\"{o}del's theorem. But, for anyone wanting
a readable, deeper and broader treatment of the theorems than I present here, I highly recommend that book.}
A short, self-contained Goldilocks exposition of a version of G\"{o}del's first theorem, aimed at a broad audience, has been lacking. 
Here I offer such an exposition.


\section{There are Non-Computable Functions}

We start with a discussion of computable and non-computable functions.
\medskip

\noindent {\bf Definition} offer We use $Q$ to denote all functions from the positive integers to $\{0,1\}$. 
That is, if $f$ is in $Q$, then for any positive integer $x$, $f(x)$ is either 0 or 1. 

Note that since
a function in $Q$ is defined on {\it all} positive integers, the number of functions in $Q$ is infinite.

\medskip

\noindent {\bf Definition} Define a function $f$ in $Q$ to be {\it computable} if there is a finite-sized computer program (in Python, for example) 
that executes on a computer (a MacBook Pro running
Snow Leopard, for example) that computes function $f$. That is, given {\it any} positive integer $x$, the program finishes in finite
time and correctly spits out the value $f(x)$. 

\medskip
\noindent {\bf Definition}
Let $A$ be the set of functions in $Q$ that are computable.

Note that the number of functions in $A$ is infinite. For example, the function $f(7) = 1$ and $f(x) = 0$ for all
$x \neq 7$ is a computable function, and we can create a similar computable function for any positive integer, in place
of 7.
So, since there are an infinite number of positive integers, there are an infinite number of computable functions.

\begin{Theorem}
\label{tnoncompute}
There are functions in $Q$ that are not computable. That is,  $A \subset Q$.
\end{Theorem}

\begin{Proof} 
In this proof we would like to talk about an {\it ordering} (or an ordered list)\footnote{The common notion of an ordering is
more technically called a {\it total order}.} of all functions in $A$,
rather than just the {\it set} $A$. 
It might seem self-evident that such an ordering should exist, and so 
one might think we could just assert that it does.
But issues of ordering are subtle; 
there are unsettled questions about which properties are sufficient to guarantee that an ordering exists.
So, we  want to be careful and fully
establish that an ordering of the functions in $A$ does exist.\footnote{I thank the students in CS 120 Fall 2014
who asked why an earlier draft of this exposition goes into such detail on the existence of an ordering.}

{\bf An ordering Exists:}
First, choose a computer language and consider a program in that language.
Each line in a program has some end-of-line symbol, so
we can concatenate the lines together into a single long string. Therefore, we think of a program in that computer
language as a single string written using some finite alphabet.

Now, since $A$ consists of the computable functions in $Q$, for each function $f \in A$, there is some computer program $P_f$
(in the chosen computer language)
that computes $f$. Program $P_f$ (considered as a single string) has some finite length. 
We can, conceptually, order the strings representing the programs that compute the functions in $A$
into a list $L$ in {\it order} of the lengths of the strings. To make the ordering perfectly precise, when there are strings of the same length, we 
order those strings
lexicographically (i.e., the way they would be {\it alphabetically} ordered in a dictionary).  
So, each {\it program} that computes a function in $A$ has a {\it well-defined} position in $L$. Then, since each function in $A$ is computed
by some program in the ordered list $L$, $L$ also defines an ordered list, which we call $L'$, containing all the functions in $A$.

A function $f$ in $A$ might be computed by different computer programs, so 
$f$ might appear in $L'$ more than once. If that occurs, we could, conceptually, 
remove all but the {\it first} occurrence of $f$ in $L'$, resulting
in an ordering of the functions in $A$, as desired.
We will see that it will not harm anything if $f$ is computed by more than one program in $L$, and hence appears in $L'$ more than once.
The only point that will matter is that there is some ordered listing $L'$ of the functions in $A$ that includes every
function in $A$.

Let $f_i$ denote the function in $A$ that appears in position $i$ in $L'$; that is, $f_i$ is computed by the $i$'th program in $L$. 
(Remember that lists $L$ and $T$ are only conceptual; we
don't actually build them--we only have to imagine them for the sake of the proof).
Next, consider a table $T$ with one column for each positive integer, and one row for each program in $L$; 
and associate the function $f_i$ with row $i$ of $T$.
Then set the value of cell $T(i,x)$ to $f_i(x)$. See Table \ref{ttableT}. 

\medskip

{\bf Function $\overline{f}$:}
Next, we define the function $\overline{f}$ from the positive integers to $\{0,1\}$ as $\overline{f}(i) = 1 - f_i(i)$.
For example, based on the functions in Table \ref{ttableT}, 
$\overline{f}(1) = 0; 
\overline{f}(2) = 1; 
\overline{f}(3) = 1; 
\overline{f}(4) = 0; 
\overline{f}(5) = 1$

Note that in the definition of $\overline{f}(i)$, 
the same integer $i$ is used both to identify
the function $f_i$ in $A$, and as the input value to $f_i$ and to $\overline{f}$. Hence the values for $\overline{f}$ are determined from
the values along the main {\it diagonal} of table $T$.
Note also that $\overline{f}$ changes 0 to 1, and changes 1 to 0. 
So, the values of function $\overline{f}$  are the {\it opposite} of the values along the main diagonal of Table $T$.
Clearly, function $\overline{f}$ is in $Q$.

\begin{table}[h!]
\begin{center}
\begin{tabular}{r|cccccccccc} 
\hline
&$1$ & $2$ & $3$ &        $4$  &                    $5$ & . & $x$ &. & $i$ & $\ldots$\\
\hline
$f_1$ & 1& 1 &        0&     0 &               0 &&$f_1(x)$&&&\\
$f_2$ & 0& 0 &        1&     0 &               0&&$f_2(x)$&&&\\
$f_3$ & 1& 1 &        0&     0 &               1&&$f_3(x)$&&&\\
$f_4$ & 0& 0 &        1&     1 &               0&&$f_4(x)$&&&\\
$f_5$ & 0& 1 &        0&     0 &               0&&$f_5(x)$&&&\\
$.$   &&&&&&.&&&&\\
$.$   &&&&&&&.&&&\\
$.$   &&&&&&&&.&&\\
$f_i$ &&&&&&&&& $f_i(i)$&\\
$\vdots$   &&&&&&&&&&$\ddots$

\end{tabular}
\caption{The conceptual Table $T$ is contains an ordered list of all computable functions in $Q$, and their values at all of the positive integers.}
\label{ttableT}
\end{center}
\end{table}

\medskip
Now we ask: Is $\overline{f}$ a computable function? 

\medskip
The answer is no for the following reason. If $\overline{f}$ were
a computable function, then there would be some row $i^*$ in $T$ such that $\overline{f}(x) = f_{i^*}(x)$ for every positive integer $x$.
For example, maybe $i^*$ is 57. But $\overline{f}(57) = 1 - f_{57}(57) \neq f_{57}(57)$, so $\overline{f}$ can't be $f_{57}$.
More generally, $\overline{f}(i^*) = 1 - f_{i^*}(i^*)$, so $\overline{f}$ and $f_{i^*}$ differ at least for one input
value (namely $i^*$), so $\overline{f} \neq f_{i^*}$. Hence, there is no row in $T$ corresponding to $\overline{f}$, and so $\overline{f}$
is not in set $A$. So $\overline{f}$ is {\it not} computable---it is in $Q$, but not in $A$.
\end{Proof}

\section{What is a Formal Proof System?}
How do we connect Theorem \ref{tnoncompute}, which is about functions, to G\"{o}del's first incompleteness theorem, which is about
logical systems?  We first must define a {\it formal proof system}.  

\medskip

\noindent {\bf Definition} A formal proof system $\Pi$ has three components: 
\begin{enumerate}
\item
A finite alphabet, 
and some finite subset words and phrases that can be used in forming (or writing) {\it statements}.\footnote{These words and phrases
are strings in the alphabet of $\Pi$. We will say that they are a subset of English.}

\item
A finite list of
{\it axioms} (statements that we take as true); and 

\item
A finite list of {\it rules of reasoning}, also called 
{\it logical inference, deduction or derivation} rules,
that can be applied to create a new statement from axioms and the statements already created, 
in an unambiguous, mechanical way. 
\end{enumerate}

The word ``mechanical" is central to the definition of rules of reasoning, and to the whole purpose 
of a formal proof system:

\begin{quote} ... we need to impose some condition to the effect that recognizing an axiom or applying a rule must
be a mechanical matter ... it is required of a formal system that in order to verify that something is an axiom or
an application of a rule of reasoning, we ... need only apply mechanical checking of the kind that can be 
carried out by a computer.\footnote{G\"{o}del's Theorem, by Torkel Franzen, CRC Press, 2005.} 
\end{quote}
\medskip

For example, the alphabet
might be the standard ASCII alphabet with 256 symbols, and Axiom 1 might be:.
``for any integer $x$, $x+1 > x$." Axiom 2 might be: ``for any integers
$x$ and $y$, $x + y$ is an integer." A derivation rule might be: ``for any three integers, $x, y, z$,
if $x > y$  and $y > z$ then $x > z$." (Call this rule the ``Transitivity Rule".) 

The finite set of allowed English words and phrases might include the phrase:
``for any integer".  Of course, there will typically be
more axioms, derivation rules, and known words and phrases than in this example.

\medskip

\subsection{What is a Formal Derivation?}

\noindent {\bf Definition} 
A {\it formal derivation} in $\Pi$ of a statement $S$ is a series of statements that begin with some 
axioms of $\Pi$, and then successively apply 
derivation rules in $\Pi$ to obtain statement $S$.

For example, $S$ might be the statement:  ``For any integer $w$,
$w + 1 + 1 > w$". A formal derivation of $S$ in $\Pi$ (using axioms and derivation rules introduced above) might be:

\begin{itemize}

\item 
$w$ is an integer, 1 is an integer, so $w+1$ is an integer (by Axiom 2).

\item
$w+1$ is an integer (by the previous statement),  1 is an integer, so $w+1 + 1 > w + 1$ (by Axiom 1).

\item
$w+1 + 1 + 1 >  w + 1 > w$ (by the previous statement and Axiom 1).

\item
$w+1 + 1 > w$ (by the Transitivity Rule). This is statement $S$.

\end{itemize}

The finite subset of English used in this formal 
derivation includes the words and phrases ``is an integer", ``by the Transitivity Rule", ``by the previous statement" etc. 
. These would be part of the finite subset of English that is part of the definition of $\Pi$.  
Each phrase used must have a clear and precise meaning in $\Pi$, so that each statement in a formal derivation,
other than an axiom, 
follows in a mechanical way from the preceding statements by the application of some derivation rule(s) or axioms.

Formal derivations are very tedious, and humans don't want to write derivations this way, but computers
can write and check them,  a fact that is key in our treatment of G\"{o}del's theorem.
(Note that what I have called a ``formal derivation" is more often
called a ``formal proof". But that is confusing, because people usually think of a ``proof" as something that establishes
a {\it true} statement, not a statement that might be false. So here we use ``formal derivation" to avoid that confusion.)

\subsection{Mechanical Generation and Checking of Formal Derivations}

We now make four key points about formal derivations.

\medskip

1. It is easy to write a program $P$ that can begin generating, 
in order of the lengths of the strings, every string $s$ that can be written in the alphabet of $\Pi$, and using
allowed words and phrases of the formal proof system $\Pi$.
Program $P$ will never stop because there is no bound on the length of the strings, and most
of the strings will not be formal derivations of anything.
But, for any finite-length string $s$ using the alphabet of $\Pi$, $P$ will eventually (and in finite time) generate $s$. 

\medskip
2. A formal derivation, being a series of statements, 
is just a string formed from the alphabet and the allowed words and phrases of the formal proof system $\Pi$. Hence, if
$s$ is any string specifying a formal derivation, $P$ will eventually (and in finite time) generate it.

\medskip
3. We can create a program $P'$ that knows the alphabet, the axioms, the deduction rules, and the meaning of the words
of the allowed subset of English used in $\Pi$, so that $P'$ 
can precisely interpret the effect of each line of a formal derivation. That is, $P'$ can {\it mechanically} check whether each line 
is an axiom, or 
follows from the previous lines by an application of some deduction rule(s) or axioms. 
Therefore, given a statement $S$, and a string $s$ that might be a formal
derivation of $S$, program $P'$ can check (in a purely mechanical way, and in finite time) 
whether string $s$ is a  formal derivation of statement $S$ in $\Pi$. 

\medskip

4. For any statement $S$, 
after program $P$ generates
a string $s$, program $P'$ can check whether $s$ is a formal derivation of statement $S$ in $\Pi$, before $P$ generates the next
string. Hence, if
there is a formal derivation $s$ in $\Pi$ of statement $S$, then $s$ will be generated and recognized in finite time
by interleaving the execution of programs $P$ and $P'$.  

\medskip

Note that  most of the strings that $P$ generates will be
garbage, and most of the strings that are not garbage will not be  formal derivations of 
$S$ in $\Pi$. But, if string $s$ is a formal derivation of statement $S$, then in finite time, program $P$ will generate $s$, 
and program $P'$
will recognize that $s$ is a formal derivation in $\Pi$ of statement $S$.

Similarly, we can have another program $P''$ that checks whether a string $s$
is a formal derivation of the statement ``not $S$", written $\neg S$. So if $\neg S$ is a statement that can be derived in $\Pi$,
the interleaved execution of programs $P$ and $P''$ will, in finite time, generate and recognize that $s$ is a formal derivation
of $\neg S$.

\section{Back to G\"{o}del}

How do we connect all this to G\"{o}del's first incompleteness theorem? 
We want to show the variant of G\"{o}del's theorem that says: in any ``rich-enough" formal proof system where no false statement about functions
can be derived, there are true statements about functions that cannot be derived.
We haven't defined what ``true" or ``rich-enough" means in general, but we will in a specific context.
\medskip

Recall function $\overline{f}$, and recall that it is well-defined, i.e., there is a value $\overline{f}(x)$ for every
positive integer $x$, and for any specific $x$, $\overline{f}(x)$ is either  0 or 1.  Recall also, that $\overline{f}$ is
not a computable function.

\medskip

\noindent {\bf Definition} We call a statement an $\overline{f}$-statement if it is either:

\begin{quote}
``$\overline{f}(x) \verb| is | 1,$"

or:

``$\overline{f}(x) \verb| is | 0,$"
\end{quote}
for some positive integer $x$. 

Note that every $\overline{f}$-statement is a statement about a specific integer. For example the statement 
``$\overline{f}(57) \verb| is | 1$" is an $\overline{f}$-statement, where $x$ has the value 57. 
Since, for any positive integer $x$, $\overline{f}(x)$ has only two possible values, 0 or 1, when the two kinds of $\overline{f}$-statements refer to the same
$x$, we refer to the first statement as $Sf(x)$ and the second statement as $\neg Sf(x)$.

\paragraph*{What is Truth?}
We say an $\overline{f}$-statement $Sf(x)$ is ``true", and $\neg Sf(x)$ is ``false",  if in fact $\overline{f}(x)$ is 1. 
Similarly, we say 
an $\overline{f}$-statement $\neg Sf(x)$  is true, and $Sf(x)$ is false, if in fact $\overline{f}(x)$ is 0. 
Clearly, for any  positive integer $x$, one of the statements $\{Sf(x), \neg Sf(x)\}$ is true and the other is false.
In this context, truth  and falsity 
are simple concepts (not so simple in general).

Clearly, it is a desirable property of a formal proof system $\Pi$, that it is not possible to give a formal derivation in
$\Pi$ for a statement that is false.

\paragraph*{What does it mean to be rich-enough?}
We need a definition.

\medskip

\noindent {\bf Definition} We define a formal proof system $\Pi$ to be {\it rich-enough} if any $\overline{f}$-statement can be formed (i.e., stated,
or written) in $\Pi$.

Note that the words ``formed", ``stated", ``written"  do not mean ``derived". The question of whether a statement can be derived in $\Pi$ is
at the heart of G\"{o}del's theorem. Here, we are only saying that the statement can be formed (or written) in $\Pi$.

\subsection{The Proof of our variant of G\"{o}del's Theorem}

Now let $\Pi$ be a rich-enough formal proof system,
and suppose {\bf a)} that $\Pi$ has the properties that no false $\overline{f}$-statements can be derived in $\Pi$; and 
suppose {\bf b)}
that for any true $\overline{f}$-statement $S$, there is a formal derivation $s$ of $S$ in $\Pi$.

Since $\Pi$ is rich-enough, for any positive integer $x$, both statements $Sf(x)$ and $\neg Sf(x)$ can be formed in $\Pi$,
and since exactly one of those statements is true, suppositions {\bf a} and {\bf b} imply that there is 
a formal derivation in $\Pi$ of exactly one of the two statements,
in particular, the statement that is true. But this leads to a contradiction of the established fact that function $\overline{f}$ is 
not computable.

In more detail, if the two suppositions ({\bf a} and {\bf b} hold, the following approach describes a computer program $P^*$ that can correctly 
determine
the value of $\overline{f}(x)$, for any positive integer $x$, in finite time.

\begin{quote}
Program $P^*$:
Given $x$, start program $P$ to successively generate all possible strings (using the finite alphabet and 
known words and phrases in $\Pi$), in order of their lengths, 
breaking ties in length lexicographically (as we did when discussing list $L$). 
After $P$ generates a string $s$, run program $P'$ to see if
$s$ is a formal derivation of statement $Sf(x)$. If it is, output that $\overline{f}(x) = 1$ and halt; 
and if it isn't, run $P''$ to see 
if $s$ is a formal derivation of $\neg Sf(x)$.  If it  is, output that $\overline{f}(x) = 0$ and halt; 
and if it isn't, let $P$ go on to generate the next possible string.
\end{quote}

The two suppositions {\bf a} and {\bf b} guarantee that for any positive integer $x$,
this mechanical computer program, $P^*$, will halt in finite time, outputting the correct value of
$\overline{f}(x)$.
But then,  $\overline{f}$ would be a {\it computable} function (computable by program $P^*$), contradicting the
already established fact that $\overline{f}$ is not a computable function.
So, the two suppositions {\bf a} and {\bf b} lead to a contradiction, so they cannot both hold for any 
rich-enough formal proof system $\Pi$.  
There are several equivalent conclusions that result. One is:

\begin{Theorem}
\label{tgodel1}
For any rich-enough formal proof system $\Pi$ in which 
no formal derivation of a false $\overline{f}$-statement is possible,
there will be some true $\overline{f}$-statement that cannot be formally derived in $\Pi$. 
\end{Theorem}

\noindent A different, but equivalent conclusion is:

\begin{Theorem}
\label{tgodel2}
In any rich-enough formal proof system $\Pi$ in which 
no formal derivation of a false $\overline{f}$-statement is possible,
there will be some positive integer $x$ such that neither statement $Sf(x)$ nor statement $\neg Sf(x)$ can be formally derived.
\end{Theorem}

We leave to the reader the proof that Theorems \ref{tgodel2} and \ref{tgodel1} are equivalent.
Theorems \ref{tgodel1} and \ref{tgodel2} are variants of G\"{o}del's first incompleteness theorem. 

\section{A Little Terminology of Formal Logic}

We proved Theorem \ref{tgodel1} with as little terminology from formal logic as possible. That was one of the goals of this exposition. Still, it
useful to introduce some terminology to be more consistent with standard use.

\medskip

\noindent {\bf Definition} A formal proof system $\Pi$ is
called {\it sound} if only true statements can be derived in $\Pi$. But note that we don't require that {\it all} true statements be derived in 
$\Pi$.

\medskip
\noindent {\bf Definition} A formal proof system $\Pi$ is called {\it complete} if for any statement $S$ that can be formed in $\Pi$,
one of  the statements $S$ or $\neg S$ can be derived in $\Pi$. But note that we don't require that the derived statement be true.

Theorem \ref{tgodel2} can then be stated as:

\begin{Theorem}
\label{tgodel3}
No formal proof system $\Pi$ that can form any $\overline{f}$ statement can be both sound and complete.
\end{Theorem}


\medskip

%
%
%
%
%

\section{G\"{o}del's Second Incompleteness Theorem}

\noindent {\bf Definition} A formal proof system $\Pi$ is called {\it consistent} if it is never possible to derive a statement $S$ in $\Pi$ and
also derive the statement $\neg S$  in $\Pi$. 

Later in the course, we will talk about G\"{o}del's second incompleteness theorem, 
which needs more machinery.  Informally, it says that if $\Pi$ is rich-enough and consistent,
there cannot be a formal derivation in $\Pi$ of the statement: ``$\Pi$ is consistent". More philosophically,
but not precisely, for any (rich enough) formal proof system $\Pi$ that is consistent, the consistency of $\Pi$ can only be established
by a different formal proof system $\Pi'$. (But then, what establishes the consistency of $\Pi'$?)

\section{Optional Homework Questions:}

1. In two places in the proofs, ties in the lengths of strings are broken lexicographically. I claim that this detail is not needed in either place.
Is this true?

\medskip

\noindent 2. In the proof of Theorem \ref{tnoncompute}, what is the point of requiring the computer programs to be listed in order of their lengths? 
Would the given proof of Theorem \ref{tnoncompute} remain correct if the programs were (somehow) listed in 
no predictable order?

\medskip

\noindent 3. In program $P^*$, what is the point of requiring program $P$ to generate strings in order of their lengths? 
Would the given proof of 
Theorem \ref{tgodel1} remain correct if $P$ did not generate the strings in that order, but could (somehow) generate all the strings in 
no predictable order?

\medskip

\noindent 4. Doesn't the following approach show that $\overline{f}(x)$ is computable?

First, create a computer program $P'$ that can look at a string $s$ over the finite alphabet used for computer programs (in some
fixed computer language, for example, C), and determine if $s$ is a legal computer program that computes a function
$f$ in $Q$. Certainly, a compiler for C can check if $s$ is a syntactically correct program in C. 

Then given any positive integer $x$,
use program $P$ to generate the strings over the finite alphabet used for computer programs, in order of their length, and
in the same order as used in table $T$. 
After each string $s$  is
generated, use program $P'$ to determine if $s$ is a program that computes a function in $Q$. Continue doing this
until $x$ such programs have been found. In terms of table $T$, that program, call it $F$, will compute function $f_x$.
Program $F$ has finite length, so $P$ will only 
generate a finite number of
strings before $F$ is generated. Then once $F$ is generated, run it with input $x$. 
By definition, program $F$ will compute $f_x(x)$ in finite time. Then output $\overline{f}(x) = 1 - f_x(x)$.

So, this approach seems to be able to compute $\overline{f}(x)$ in finite time, for any positive integer $x$, 
showing that $\overline{f}$ is a computable function.
Doesn't it?

Discuss and resolve.

\medskip

\noindent 5. Use the resolution to the issue in problem 4, to state and prove an interesting theorem about computer programs (yes, this is
a vague question, but the kind that real researchers face daily).

\medskip

\noindent 6. Show that theorems \ref{tgodel1} and \ref{tgodel2} are equivalent.

\medskip

\noindent 7. Show that a formal proof-system that is sound is also consistent. Then ponder whether it is true that any formal proof-system that
is consistent must be sound. Hint: no.  


\section{What is not in this exposition?}

Lots of stuff that you might see in other proofs and expositions of G\"{o}del's first incompleteness theorem:
propositional and predicate logic, models, WFFs, prime numbers, prime factorization theorem, Chinese remainders,
G\"{o}del numbering,  countable and uncountable infinity,
self-reference, recursion, paradoxes, liars, barbers, librarians, ``This statement is false", 
Peano postulates, Zermelo-Fraenkel set theory, Hilbert, Russell, Turing, universal Turing machines, the halting problem, undecidability, ..., 
quantum theory, insanity, neuroscience, the mind, zen,  self-consciousness, evolution, relativity, philosophy, religion, God, 
Stalin, ... . The first group
of topics are actually related in precise, technical ways to the theorem, but can be avoided, as done in this exposition.\footnote{Most
expositions of G\"{o}del's theorems use self-reference, which I find unnecessarily head-spinning, and I think its use is sometimes
intended to make G\"{o}del's theorem seem deeper and more mystical than it already is.}
Some of those related technical topics are important in their own right, particularly Turing undecidability, which we will cover
in detail later in the course. 
The second group of topics are not related in a precise, technical way to the theorem. Some are fascinating in their own
right, but their inclusion makes G\"{o}del's theorem more mystical, and should not be confused for its actual technical content.

\section{Final Comments}
The exposition here does not follow G\"{o}del's original proof, and while the exposition is my own, the general
approach reflects (more and less) 
the contemporary computer-sciencey way that G\"{o}del's theorem is thought about, i.e., via computability.  In coming to this exposition for undergraduates, 
I must acknowledge the discussion 
of G\"{o}del's theorems in Scott Aaronson's book {\it Quantum Computing Since Democritus}, and an exposition shown to me by Christos
Papadimitriou. Those are both shorter, aimed at a more advanced audience, and are based on the undecidability of Turing's Halting problem. I also
thank David Doty for pointing out an incorrect definition in my first version of this exposition.

Second, I must state again that the variant of G\"{o}del's theorem proved here is weaker
than what G\"{o}del originally proved. In this exposition, the formal proof system must be able to express any 
$\overline{f}$-statement,
but G\"{o}del's original proof only requires that the formal system 
be able to express statements about arithmetic (in fact, statements about arithmetic on integers only using 
addition and multiplication).  This 
is a more limited domain, implying
a stronger theorem. That difference partly explains why a proof of G\"{o}del's original theorem is technically more demanding than the
exposition here. Further, G\"{o}del did not just prove the {\it existence} of a true statement that could not be derived (in any sufficiently rich, 
sound proof system), he demonstrated a {\it particular} statement with that property. But, I believe that the moral, cultural, mathematical,
and philosophical impact of the variant of G\"{o}dels theorem proved here 
is comparable to that of G\"{o}del's actual first incompleteness theorem. Many modern treatments of G\"{o}del's theorem similarly reflect this view. 
Of course, some people disagree and insist that anything using the phrase ``G\"{o}del's theorem" must
actually be the same as what G\"{o}del proved.\footnote{Although, G\"{o}del did not actually prove what is generally stated as 
``G\"{o}del's theorem", but only a weaker form of it--which was later strengthened by Rosser to become the classic
``G\"{o}del's theorem".  Accordingly, some people call it the  ``G\"{o}del-Rosser theorem".}
\end{document}